\documentclass[10pt,reqno]{amsart}

\usepackage{amssymb}
\usepackage{amsfonts}
\usepackage{amsthm}
\usepackage{amsmath}
\usepackage{bbm}
\usepackage{xcolor}
\usepackage{enumerate}
\usepackage{mathabx}

\numberwithin{equation}{section}
\allowdisplaybreaks

\newtheorem{theorem}{Theorem}[section]
\newtheorem{proposition}[theorem]{Proposition}

\newtheorem{lemma}[theorem]{Lemma}
\newtheorem{corollary}[theorem]{Corollary}

\theoremstyle{definition}
\newtheorem{definition}[theorem]{Definition}

\theoremstyle{remark}
\newtheorem{remark}[theorem]{Remark}

\newcommand{\R}{\mathbb{R}}
\newcommand{\Rd}{\mathbb{R}^d}

\newcommand{\N}{\mathbb{N}}
\newcommand{\Z}{\mathbb{Z}}

\renewcommand{\hat}{\widehat}
\newcommand{\eps}{\varepsilon}

\newcommand{\scriptD}{\mathcal{D}}
\newcommand{\scriptE}{\mathcal{E}}

\newcommand{\scriptQ}{\mathcal{Q}}
\newcommand{\scriptR}{\mathcal{R}}
\newcommand{\scriptS}{\mathcal{S}}

\DeclareRobustCommand{\rchi}{{\mathpalette\irchi\relax}}
\newcommand{\irchi}[2]{\raisebox{\depth}{$#1\chi$}} 

\newcommand{\ceil}[1]{\left\lceil #1 \right\rceil}

\newcommand{\qtq}[1]{\quad\text{#1}\quad}

\DeclareMathOperator*{\supp}{supp}

\DeclareMathOperator*{\dist}{dist}
\DeclareMathOperator*{\wklim}{wk-lim}

\begin{document}
\title[Existence of extremizers for Fourier restriction to the moment curve]{Existence of extremizers for Fourier restriction to the moment curve}
\author{Chandan Biswas and Betsy Stovall}
\address{Department of Mathematics, Indian Institute of Science, Bangalore, 560012}
\email{chandanbiswa@iisc.ac.in}
\address{Department of Mathematics, University of Wisconsin, Madison, WI 53706}
\email{stovall@math.wisc.edu}

\begin{abstract}
We show that the restriction and extension operators associated to the moment curve possess extremizers and that $L^p$-normalized extremizing sequences of these operators are precompact modulo symmetries.  
\end{abstract}

\maketitle


\section{Introduction}


This article establishes the existence of extremizers and compactness, modulo symmetries, for the restriction/extension inequalities associated to the moment curve. More precisely, we consider the operator 
$$
\scriptE f(x) := \int_\R e^{i x \cdot \gamma(t)} f(t)\, dt,  \qquad \gamma(t) := (t, t^2, \ldots, t^d),
$$
which was shown by Drury~\cite{Drury85} to extend as a bounded linear operator from $L^p(\R)$ to $L^q(\R^d)$ if and only if $q > \tfrac{d^2 + d + 2}2$ and $q = \tfrac{d (d + 1)}2 p'$ for $d \geq 2$. We prove that for all $(p, q)$ in this range, there exist nonzero functions $f$ such that $\|\scriptE f \|_q = \|\scriptE\|_{L^p \to L^q} \|f\|_p$. Moreover, whenever $(p, q) \neq (1, \infty)$, $L^p$-normalized extremizing sequences (i.e. those that saturate the operator norm) possess subsequences that converge, modulo the application of symmetries of the operator, to an extremizing function.

Our argument uses a modified version of the concentration-compactness framework of Lions~\cite{Lions85} and the related Method of Missing Mass of Lieb~\cite{Lieb83}. Such methods have been well-studied for the $L^2$-based restriction/extension problems associated to certain hypersurfaces (see \cite{Carneiro2008}, \cite{CarneiroFoschiSilvaThiele2015}, \cite{CarneiroDiogo2014}, \cite{ChristShao2010}, \cite{Foschi2004}, \cite{frankLiebSabin2016}), and the resulting theory has been an important step towards breakthroughs in the study of long-time behavior of various dispersive equations, including NLS, NLW, and other equations (\cite{BrocchiSilvaQuilodran2020}, see also \cite{FoschiSilva2017}, \cite{KillipVisan2013} and the references therein).   

In this article, we make two advances relative to these previous works.  First, there are not, to our knowledge, any previous results in the literature regarding concentration-compactness phenomena for Fourier restriction to higher co-dimensional manifolds.  Though Fourier restriction to lower-dimensional manifolds has been less intensively studied than restriction to hypersurfaces, one is naturally led to the former from the latter by examining the sublevel sets of the Gaussian curvature of certain higher-order surfaces.  (Some preliminary work on this connection is in~\cite{SchwendStovall19}.)  Naturally, we begin this study by examining the model ``curved curve,'' which is referred to as the \textit{moment curve} in the literature.  A key step is a multilinear generalization of the bilinear-to-linear argument of Tao--Vargas--Vega \cite{TVV} (and the later refinement thereof by B\'egout--Vargas~\cite{BV}), for which a new Whitney-like decomposition is needed.  Second, we continue the development from \cite{Stovall20} of effective concentration compactness techniques for $L^p \to L^q$ inequalities for general exponent pairs $(p,q)$.  While \cite{Stovall20} laid out a strategy that makes accessible exponent pairs with $p \neq 2$, a new issue arises in the higher co-dimensional setting.  Namely, the high degree of the determinant in the scaling relation $\gamma(\lambda t) = \mathop{diag}(\lambda,\ldots,\lambda^d)\gamma(t)$ puts certain exponent pairs $(p,q)$ outside of the range accessible by direct adaptations of the Tao--Vargas--Vega approach.  We bridge the gap between this multilinear-accessible range ($1 < p < d + 2$) and the Drury range ($1 < p < \frac{d^2 + d + 2}2$) by introducing a sort of interpolation-like argument, which seems likely to be of use in other settings.  

We now turn to the precise formulation of our results, for which we introduce some notation and terminology.  By a \textit{symmetry} of the operator $\scriptE : L^p \to L^q$, we mean an element $S$ of the isometry group of $L^p(\R)$ for which there exists a corresponding element $T$ of the isometry group of $L^q(\R^d)$ obeying $\scriptE \circ S = T \circ \scriptE$. The key symmetries for our analysis are the dilations, the frequency translations, and the modulations. More precisely, the dilations are given by: 
$$
f \mapsto f^\lambda := \lambda^{- \frac1p} f(\lambda^{-1} \cdot), \quad \scriptE f^\lambda(x) = \lambda^{\frac{d (d + 1)}{2q}} \scriptE f(D_\lambda(x)), \quad D_\lambda(x) = (\lambda x_1, \ldots, \lambda^d x_d);
$$
The translations are given by:
\begin{equation}\label{define A_t}
f \mapsto \tau_{t_0} f := f(\cdot - t_0), \qquad \scriptE(\tau_{t_0} f)(x) = \mathbb{L}_{t_0} (\scriptE f)(x),
\end{equation}
where $\mathbb{L}_{t_0}$ is the \textit{boost} $\mathbb{L}_{t_0} g(x) := e^{i x \cdot \gamma(t_0)} g(A_{t_0}^T x)$, and $A_{t_0}$ is the unique element in $GL(d)$ for which $\gamma(t + t_0) = \gamma(t_0) + A_{t_0} \gamma(t)$ (we observe that $A_{t_0}$ is lower triangular, with ones on the diagonal); The modulations are given by:
$$
(m_{x_0} f)(t) := e^{- i x_0 \cdot \gamma(t)} f(t), \qquad \scriptE(m_{x_0} f) = \tau_{x_0} \scriptE f.
$$

Let $1 \leq p < \tfrac{d^2 + d + 2}2$, $q:=\tfrac{d(d+2)}2 p'$, and set $B_p := \|\scriptE\|_{L^p \to L^q}$. We say that $f \in L^p$ is an extremizer of (the $L^p(\R) \to L^q(\R^d)$ inequality for) $\scriptE$ if $f \not \equiv 0$ and $\|\scriptE f\|_{q} = B_p \|f\|_p$. We say that $\{f_n\} \subseteq L^p$ is an extremizing sequence for (the $L^p \to L^q$ inequality for) $\scriptE$ if $f_n \not \equiv 0$ for all $n$ and $\lim_{n \to \infty} \tfrac{\|\scriptE f_n\|_q}{\|f_n\|_p} = B_p$. We are most interested in normalized extremizing sequences, that is, extremizing sequences $\{f_n\}$ with $\|f_n\|_p = 1$, for all $n$. Our main result is the following.

\begin{theorem}\label{mainthm}
For $d \geq 2$ there exist extremizers of the $L^p(\R) \to L^q(\R^d)$ inequality for $\scriptE$ for every $(p, q)$ satisfying $1 \leq p < \frac{d^2 + d + 2}{2}$ and $q = \tfrac{d (d + 1)}2 p'$. Moreover, when $p > 1$, given any extremizing sequence of $\scriptE$, there exists a subsequence that converges to an extremizer in $L^p(\R)$, after the application of a suitable sequence of symmetries.  
\end{theorem}

We note that the existence of extremizers and non-compactness of extremizing sequences in the case $p=1$ is elementary.  Theorem~\ref{mainthm} immediately yields a related result for the corresponding restriction operator.  

\begin{corollary}
The analogous result holds for the restriction operator $\scriptR g(t) := \hat g (\gamma(t)))$; namely, for $d \geq 2$ there exist extremizers of the $L^r(\R^d) \to L^s(\R)$ inequality for $\scriptR$ for every $(r, s)$ satisfying $1 \leq r < 1 + \frac 2{d^2 + d}$ and $r' = \tfrac{d^2 + d}2 s$. Moreover, when $r > 1$, given any extremizing sequence of $\scriptR$, there exists a subsequence that, after the application of a suitable sequence of symmetries, converges in $L^r(\R)$ to an extremizer.  
\end{corollary}

We give the short proof of the corollary now.  
\begin{proof}
In the case $r = 1$, the first conclusion is elementary, and so we assume that $r > 1$. Let $\{g_n\}$ be an $L^r$-normalized extremizing sequence of $\scriptR$. Let $f_n := B_{s'}^{- (s - 1)} |\scriptR g_n|^{s - 2} \scriptR g_n$. Then $\lim_n \|f_n\|_{s'} = 1$. On the other hand, by duality,
$$
B_{s'} = B_{s'}^{- (s - 1)} \lim_n \|\scriptR g_n\|_s^s = \lim_n  \langle g_n, \scriptE f_n \rangle \leq \|g_n\|_r \|\scriptE f_n\|_{r'} \leq B_{s'},
$$
so $\{f_n\}$ is an extremizing sequence for $\scriptE$, the norms of whose elements tend to $1$. Applying Theorem~\ref{mainthm} to $\{\|f_n\|_{s'}^{- 1} f_n\}$, there exist symmetries $\{S_n\}$ of $\scriptE$ such that along a subsequence, $\{S_n f_n\}$ converges in $L^{s'}$ to some extremizer $f$ for $\scriptE$. Let $T_n$ denote the corresponding $L^{r'}$ automorphism, that is, $T_n \circ \scriptE = \scriptE \circ S_n$. Since $S_n f_n = B_{s'}^{- (s - 1)} |\scriptR T_n g_n|^{s - 2}\scriptR T_n g_n$, replacing $\{g_n\}$ by $\{T_n g_n\}$ if necessary, we may assume that $S_n$ equals the identity for all $n$ and that $\{f_n\}$ converges to $f$. By Banach-Alaoglu, $\{g_n\}$ converges weakly to some $g \in L^r$ and so $\|g\|_r \leq \lim_n \|g_n\|_r = 1$. On the other hand
$$
B_{s'} = \lim_n \langle g_n,\scriptE f_n \rangle = \langle g, \scriptE f \rangle \leq B_{s'} \|g\|_r,
$$
so $\|g\|_r = 1$. By Theorem $2.11$ in \cite{LL} $\{g_n\}$ converges in $L^r$ to $g$, from which we additionally conclude that $g$ is an extremizer of $\scriptR$.
\end{proof}

\subsection*{Acknowledgements}  While conducting this research, the first named author was supported by C.~V.\ Raman Postdoctoral fellowship, and the second named author was partially supported by NSF grant DMS-1653264 and the Wisconsin Alumni Research Foundation (WARF).  The authors are indebted to Benjamin Bruce and the anonymous referee for their comments on earlier versions of this manuscript.  

\section{Outline of proof}

We follow the general outline laid out in \cite{Stovall20}, which consists of two key steps: First, we show that an extremizing sequence possesses good frequency (i.e.\ along  $\R$) localization, after application of a suitable sequence of symmetries (translations and dilations). Second, we show that extensions of an extremizing sequence have good spatial localization (after modulation). This enables us to upgrade weak convergence to $L^p$ convergence using the fact that our sequence is extremizing.  The $L^p$ limit is necessarily an extremizer.  

For the frequency localization, our first step is to prove that a nonnegligible contribution to each $\scriptE f_n$ comes from a single well-localized piece, $f_n \rchi_{I} \rchi_{|f_n|^p |I| \lesssim 1}$, stated more precisely as Proposition~\ref{P:refined Drury}. To prove this, we develop a $d$-linear-to-linear version of the bilinear-to-linear argument of Tao--Vargas--Vega (\cite{TVV}), with improved efficiency in the spirit of Begout--Vargas \cite{BV}.  Though strong $d$-linear adjoint restriction theorems predate (and are integral in the proof of) Drury's Theorem, the higher order linearity presents some new geometric challenges (relative to the bilinear case) as we implement them to detect well-localized pieces of the $f_n$.  In particular, developing a Whitney decomposition of the off-diagonal in $\R^d$ requires determining a notion of relatedness for $d$-tuples of intervals (as opposed to pairs of balls, which arise in the bilinear setting).   Furthermore, for basic arithmetic reasons relating to the magnitude of the scaling factor $\tfrac{d(d+1)}2$ in the relation $q=\tfrac{d(d+1)}2 p'$ (compared with $\tfrac{d+2}d$ in the case of elliptic hypersurfaces), a straightforward adaptation of the methods of Tao--Vargas--Vega cannot yield scale-invariant $L^p \to L^q$ inequalities for large values of $p$ ($p>d+2$).  We circumvent this difficulty by adapting the proof of the Marcinkiewicz interpolation theorem.  

Having identified a single contributor to much of $\scriptE f_n$, iteration yields a bounded number of contributors to any specified proportion of $\scriptE f_n$ (Lemma~\ref{L:nibble}).  However, localization requires a bit more, namely, that (after applying symmetries) these pieces must all be at scale and of maximum magnitude about one (Proposition~\ref{P:freq loc}).  Due to convexity (i.e., $q>p$), this follows by proving an orthogonality result (Lemma~\ref{L:Lq freq orthog}), utilizing that ``distant'' pieces (i.e., those with disparate localizations) have extensions that interact weakly. We use bilinear estimates based on either H\"older's inequality or a multilinear inequality of Christ from \cite{christ_curve} to establish weak interactions.  

Finally, having established that for any extremizing sequence $\{f_n\}$ is (after applying symmetries), $\{|f_n|\}$ is well-approximated in a uniform way by uniformly bounded, compactly supported functions, to prove convergence, we need to control the oscillations of the $f_n$.  To this end, we develop a profile decomposition result (Proposition~\ref{P:profile}) to write the bounded, compactly supported approximations as a superposition of modulated profiles (the modulations, but not the profiles, may depend on $n$).  Curvature and stationary phase enable us to prove an $L^p$ almost orthogonality result for these profiles, and so for $\{f_n\}$ extremizing, there is exactly one significant profile.  Using basic properties of $L^p$ spaces, we can then remove the truncations of the $f_n$ (to bounded, compactly supported functions), without disturbing our profiles nor modulations too much.  We thus obtain $L^p$ convergence of the $f_n$ to an extremizer, our desired outcome.

\subsection*{Future directions} We believe that many of our methods have the potential to be extended to a larger class of curves (and indeed, the authors intend to do so in a forthcoming article).  However, our proof uses the symmetries for the moment curve in a fundamental way in the passage from the multilinear to the linear inequality (the Whitney decomposition step, in particular), and a number of changes would be needed to extend this to more general curves, even when the torsion is comparable to 1.  Therefore such results are outside of the scope of this article.  Finally, we note that the analogous questions for manifolds of intermediate dimension (dimension and co-dimension both strictly larger than 1) seem to be extremely interesting.

\subsection*{Notation} We write $A \lesssim B$ to denote $A \leq C B$ where $C$ may depend on the dimension $d$ and the exponent $p$, and whose value may change from one line to the next but is independent of $A$ and $B$. For the rest of the article we assume that $d \geq 3$.  

\section{A refined extension estimate} \label{S:refined extn}

The purpose of this section is to prove two refinements (Propositions~\ref{P:refined Drury} and~\ref{P:refined Drury 2}) of the $L^p \to L^q$ inequalities of Drury, both of which will be used in the proof of Theorem~\ref{mainthm}.   These results show that if $f$ has nonnegligible extension, then a significant portion of the extension comes from a piece of $f$ with good frequency localization.  Later, we will capture essentially all of the extension of $f$ by iterating this inequality.  

\begin{proposition} \label{P:refined Drury}
Let $1 < p < \tfrac{d^2 + d + 2}2$, and let $q := \tfrac{d (d + 1)}2 p'$. There exist $0 < \theta = \theta_p < 1$ and $c_p > 0$ such that 
\begin{equation} \label{E:refined Drury}
\|\scriptE f\|_q \lesssim \left(\sup_{k \in \Z} \sup_{I \in \scriptD_k} \sup_{n \geq 0} 2^{- c_p n} \|f_I^n\|_p \right)^{1 - \theta} \|f\|_p^{\theta}, \qquad f \in L^p.
\end{equation}
Here $\scriptD_k$ denotes the set of all dyadic intervals of length $2^k$, and 
$$
f_I^n := f \rchi_{I} \rchi_{\{|f| < 2^n \|f\|_p |I|^{- \frac1p}\}}.  
$$
\end{proposition}

\begin{proof}[Proof of Proposition~\ref{P:refined Drury}]
We will prove the proposition in two steps : First, in the range $p < d + 2$, we will prove the inequality (\ref{E:refined Drury}) by using a multilinear extension estimate and a variant of an argument of B\'egout--Vargas \cite{BV}; then, we will adapt real interpolation methods to deduce this bound for larger values of $p$. The significance of $p < d + 2$ is that it ensures that $(\frac qd)' > \frac pd$, which allows for a $d$-linear-to-linear variant of the bilinear-to-linear argument of Tao--Vargas--Vega~\cite{TVV}. We start with the following lemma.

\begin{lemma}\label{L:multilinear extn}
Let $I_1, \ldots, I_d$ be intervals of length one, and assume that there exists some $k$, $1 \leq k < d$, such that for all $j \leq k$ and $j' > k$, $\dist(I_j, I_{j'}) \gtrsim 1$. Then for $f_j$ supported on $I_j$ and $q > \tfrac{d^2 + d + 2}2$,
\begin{equation} \label{E:multilinear extn}
\|\prod_{j = 1}^d \scriptE f_j\|_{\frac qd} \lesssim \prod_{j = 1}^d \|f_j\|_s, \qquad s := (\tfrac{2 q}{d^2})'.
\end{equation}
\end{lemma}

\begin{proof}
Changing variables,
$$
\prod_{j = 1}^d \scriptE f_j(x) = \int F(\xi) e^{i \xi \cdot x}\,d\xi,
$$
where $\xi = \sum_{j = 1}^d \gamma(t_j)$, $t_1 < \cdots < t_d$, $F(\xi) = \sum_{\sigma \in S_d} \prod_{j = 1}^d f_{\sigma(j)}(t_j) \prod_{i < j\leq d} |t_i - t_j|^{-1}$, and $S_d$ denotes the symmetric group on $d$ letters.

We set $a := (\frac qd)'$. Since $a < 2$ we may apply Hausdorff--Young to see that
\begin{align}\notag
\big(\|\prod_{j = 1}^d \scriptE f_j\|_{\frac qd}\big)^a
&\lesssim \|F\|^a_a
= \int \prod_{j = 1}^d |f_j|^a(t_j) \prod_{i < j \leq d} |t_i - t_j|^{- (a - 1)} dt
\\\label{E:break}
&\lesssim \Big(\int \prod_{j = 1}^k |f_j|^a(t_j) \prod_{i < j \leq k} |t_i - t_j|^{- (a - 1)} dt_1 \cdots dt_k\Big)
\\\notag
&\qquad \times \Big(\int\prod_{j = k + 1}^d |f_j|^a(t_j) \prod_{k + 1\leq i < j \leq d} |t_i - t_j|^{- (a - 1)} dt_{k + 1} \cdots dt_d\Big).
\end{align}

Our conditions on $q$ imply that $a - 1 = \tfrac{a}{a'} < \tfrac2{d - 1}$. Thus by Proposition~2.2 of \cite{christ_curve}, the right hand side of \eqref{E:break} is bounded by
$$
\bigl(\prod_{j = 1}^k \|f_j\|_{b_ka}^a \bigr)\bigl(\prod_{j = k + 1}^d \|f_j\|_{b_{d - k}a}^a \bigr), \qquad b_n := \bigl(\tfrac2{(n - 1)(a - 1)}\bigr)'.
$$
After a bit of arithmetic, we see that $b_n a \leq \bigl(\tfrac{2 q}{d^2})'$, whenever $1 \leq n < d$, so \eqref{E:multilinear extn} follows from H\"older's inequality, since each $f_j$ is supported on a set of measure at most one.
\end{proof}

By scaling, Lemma~\ref{L:multilinear extn} immediately implies the following corollary.

\begin{lemma}\label{L:scaled multilinear extn}
Let $q > \frac{d^2 + d + 2}2$ and let $I_1, \ldots, I_d$ be intervals of length $r > 0$, and assume that there exists some $k$, $1 \leq k < d$, such that for $j \leq k < j'$, $\dist(I_j, I_{j'}) \gtrsim r$. Then for functions $f_j$ supported on $I_j$, $1 \leq j \leq d$,
\begin{equation} \label{E:scaled multilinear extn}
\|\prod_{j = 1}^d \scriptE f_j\|_{\frac qd} \lesssim r^{d(\frac1{s'} - \frac{d (d + 1)}{2 q})} \prod_{j = 1}^d \|f_j\|_s, \qquad s :={(\tfrac{2 q}{d^2})'}.
\end{equation}
\end{lemma} \qed

In the spirit of the bilinear-to-linear argument of Tao--Vargas--Vega \cite{TVV}, we turn to a Whitney decomposition of $\R^d$ on whose pieces we can apply bound~\eqref{E:scaled multilinear extn}.  

Consider the diagonal $\Delta := \{(t, \ldots, t) : t \in \R\}$ and the annular tubes $T_r := \{\xi: \tfrac12 r \leq \dist(\xi, \Delta) \leq 2 r\}$. We cover $T_{2^m}$ with axis-parallel dyadic cubes of side length $2^{m - K_d}$, with $K_d$ sufficiently large for later purposes. Let $\scriptQ$ denote a maximal nonoverlapping collection of such dyadic cubes, and let $\scriptQ_n$ denote the subcollection consisting of those cubes in $\scriptQ$ having sidelength $2^n$.  We may assume that the collection $\scriptQ$ is invariant under permutations of the coordinates. Each $Q \in \scriptQ_n$ may be written
\begin{equation} \label{E:Q via symmetries}
Q = I^Q_1 \times \cdots \times I^Q_d = 2^n(Q_{\vec l} + \vec k), \qtq{with} Q_{\vec l} := [0,1]^d + \vec l,
\end{equation}
where the $I^Q_j$ are intervals, $\vec k, \vec l \in \Z^d$, and
$$
\vec k = (k,\ldots,k), \qtq{and} \vec l \in [0,N_d]^d \setminus ([0,M_d]^d + \Delta),
$$
with $M_d < N_d$ large dimensional constants depending on $K_d$.  The expression \eqref{E:Q via symmetries} is uniquely determined by $Q$ if we require that some entry of $\vec l$ equals 0.    

We note that
\begin{equation} \label{E:Ef^d expansion}
\|\scriptE f\|_q^d = \|(\scriptE f)^d\|_{q/d} = \|\sum_{Q \in \scriptQ} \prod_{j=1}^d \scriptE (f\rchi_{I^Q_j})\|_{q/d}.  
\end{equation}
After a possible (harmless) reordering of indices and pigeonholing, we see that the hypotheses of Lemma~\ref{L:scaled multilinear extn} apply to each $Q_{\vec l}$, and thus to each $Q \in \scriptQ$.  However, to access the summands on the right hand side of \eqref{E:Ef^d expansion} for application of Lemma~\ref{L:scaled multilinear extn}, we need a bit more.  We turn now to an adaptation of the Whitney decomposition and almost orthogonality argument of \cite{TVV}.  

Define $\Gamma(t_1, \ldots, t_d) := \sum_{j = 1}^d \gamma(t_j)$. Then for $Q \in \scriptQ$, $\prod_{j=1}^d \scriptE (f\rchi_{I^Q_j})$ has Fourier support contained in $\Gamma(Q)$.  With $Q$ as in \eqref{E:Q via symmetries}, with some entry of $\vec l$ equal to 0, we observe that
\begin{equation} \label{E:GammaQ via symmetries}
\Gamma(Q) = D_{2^n}(A_k \Gamma(Q_{\vec l}) + d \gamma(k)),
\end{equation}
where $A_k$ is from~\eqref{define A_t}. This motivates us to define a map 
$$
S_Q(\zeta):=D_{2^n}(A_k \zeta + d\gamma(k)).
$$  
We will use the $S_Q$ in a fundamental way in the proof of the decomposition lemma below.  

\begin{lemma} \label{L:Whitney}
There exists a collection $\{\psi_Q\}_{Q \in \scriptQ}$ of smooth functions with the following properties:   $\sum_Q \psi_Q \equiv 1$ on $\Gamma(\R^d \setminus \Delta)$, the support of each $\psi_Q$ intersects the support of a bounded number of other $\psi_{Q'}$, each $\Gamma(Q)$ intersects the support of a bounded number of supports of $\psi_Q$, and $\{\|\widecheck\psi_Q\|_{L^1}\}$ is a bounded set.
\end{lemma}

\begin{proof}[Proof of Lemma~\ref{L:Whitney}]
For each $\vec l \in \Z^d \cap [0,N_d]^d \setminus ([0,M_d]^d + \Delta)$, we let $V_{\vec l}$ denote a neighborhood of $\Gamma(Q_{\vec l})$, sufficiently small for later purposes.  For $Q$ taking the form  \eqref{E:Q via symmetries}, with some entry of $\vec l$ equal to 0, we define $V_Q:=S_Q(V_{\vec l})$, a neighborhood of $\Gamma(Q)$.  We claim that the $V_Q$ are finitely overlapping.  More precisely, we will show that if $Q$ takes the form \eqref{E:Q via symmetries} and $Q' = 2^{n'}(Q_{\vec l'} + \vec k')$, with $\vec k', \vec l'$ satisfying conditions analogous with $\vec k, \vec l$, then $V_Q \cap V_{Q'} \neq \emptyset$ implies $|n-n'| \lesssim 1$ and $|k-2^{n'-n}k'| \lesssim 1$.  

To this end, we define 
$$
\rho(\xi) := \sum_{i = 1}^d |\xi_i|^{\frac1i}, \qquad \delta(\xi) := \min_{t \in \R} \rho(A_{- t} (\xi - d \gamma(t))),
$$
and let $t(\xi)$ denote the minimum of all $t$ with $\delta(\xi) = \rho(A_{- t} (\xi - d \gamma(t)))$. (By basic calculus, we can see that these minima are attained.)  We observe that if $\xi \in \Gamma(Q_{\vec l})$, then $|t(\xi)| \lesssim 1$ and $\delta(\xi) \sim 1$. This is because $\Gamma(Q_{\vec l})$ is compact and does not intersect $d \cdot \gamma$, provided $K_d$ is sufficiently large. Therefore, we may choose $V_{\vec l}$ sufficiently small so that $|t(\xi)| \lesssim 1$ and $\delta(\xi) \sim 1$ for every $\xi \in V_{\vec l}$.  If $\xi \in V_Q = S_Q(V_{\vec l})$, then (after some basic linear algebra) $\xi = S_Q\zeta = A_{2^{n}k}D_{2^n} \zeta + d\gamma(2^{n}k)$, for some $\zeta \in V_{\vec l}$, so 
$$
\rho(\xi) = \rho(D_{2^n}\zeta) = 2^n \rho(\zeta) \sim 2^{n}, \qtq{and} t(\xi) = t(D_{2^n}\zeta) + 2^n k = 2^n(t(\zeta)+k),
$$
and the latter implies that $|2^{-n}t(\xi)-k| \lesssim 1$.  If, in addition, $\xi \in V_{Q'}$, the same computations and an application of the triangle inequality imply our claim that $2^{n'-n} \sim 1$ and $|2^{n'-n}k'-k| \lesssim 1$.  

We now determine our $\psi_Q$.  For $\vec l \in \Z^d \cap [0,N_d]^d \setminus ([0,M_d]^d + \Delta)$, let $\phi_{\vec l}$ denote a smooth, nonnegative function, identically 1 on $Q_{\vec l}$ and identically 0 off of $V_{\vec l}$.  For $Q = S_Q(V_{\vec l})$, we set $\phi_Q:=\phi_{\vec l} \circ (S_Q)^{-1}$.  Then $\sum_Q \phi_Q \sim 1$, and the sum has a bounded number of nonzero entries at each point.  We define 
$$
\psi_Q :=\phi_Q(\sum_{Q' : V_Q \cap V_{Q'} \neq \emptyset}\phi_{Q'})^{-1}.  
$$
The support and partition of unity conditions from the lemma are immediate.  For the $L^1$ bound on the $\widecheck\psi_Q$, we note that from the computations above, the set $\{S_Q^{-1}S_{Q'} : V_Q \cap V_{Q'} \neq \emptyset\}$ is precompact in the set of invertible affine transformations, whence the set $\{\phi_Q \circ S_Q':V_Q \cap V_{Q'} \neq \emptyset\}$ is precompact in the Schwartz class.  Therefore $\{(\psi_Q \circ S_Q)\widecheck{\:}\}$ is precompact in $\scriptS$ and consequently bounded in $L^1$.  
\end{proof}

Using Lemma~\ref{L:scaled multilinear extn} and almost orthogonality, we will prove the following.   

\begin{lemma} \label{L:Xp,dt,dt,s}
If $q > \tfrac{d^2 + d + 2}2$, then
$$
\|\scriptE f\|_q^d \lesssim \bigl(\sum_n \sum_{I \in \scriptD_n} 2^{n t d (\frac1{s'} - \frac{d (d + 1)}{2 q})} \|f_I\|_s^{d t} \bigr)^{\frac1t},
$$
where $s := (\frac{2 q}{d^2})'$, $t :=(\frac qd)'$, $\scriptD_n$ is the collection of all dyadic intervals of length $2^n$, and $f_I := f \rchi_I$.  
\end{lemma}

\begin{proof}
We first state the almost orthogonality result that we need, a slight modification of \cite[Lemma 6.1]{TVV} with the same proof.  By interpolating the cases $p=1,2,\infty$, where the result is elementary (triangle inequality and Hausdorff--Young for $p=1,\infty$, Plancherel for $p=2$), the operator $T\Big(\{g_Q\}_{Q \in \scriptQ}\Big):=\sum_{Q \in \scriptQ} g_Q * \widecheck\psi_Q$ maps $\ell^{\tilde p'}(L^p)$ boundedly into $L^p$, $1 \leq p \leq \infty$, where $\tilde p':=\min\{p,p'\}$.  Noting that $\widetilde{q/d} = (q/d)'=t$, we have by \eqref{E:Ef^d expansion}, $\sum \psi_Q \equiv 1$, and the finite overlap condition
\begin{align*}
\|\scriptE f\|_q^d 
&= \|\sum_Q \sum_{Q':V_Q \cap V_{Q'} \neq \emptyset}  \widecheck\psi_Q * \prod_{j=1}^d \scriptE (f\rchi_{I^{Q'}_j})\|_{q/d} \\
&\lesssim \bigl(\sum_Q \|\sum_{Q':V_Q \cap V_{Q'} \neq \emptyset}  \prod_{j=1}^d \scriptE (f\rchi_{I^{Q'}_j})\|_{q/d}^t\bigr)^{1/t}
 \lesssim \bigl(\sum_n \sum_{Q \in \scriptQ_n} \|\prod_{j = 1}^d \scriptE f_{I_j} \|_{\frac qd}^t \bigr)^{\frac 1t}.  
\end{align*}
Applying Lemma~\ref{L:scaled multilinear extn}, we obtain
\begin{align}\label{E:Xp,dt,dt,s}
\|\scriptE f\|_q^d 
\lesssim \bigl(\sum_n \sum_{Q \in \scriptQ_n} 2^{n t d (\frac1{s'} - \frac{d (d + 1)}{2 q})} \prod_{j = 1}^d \|f_{I_j}\|_s^t \bigr)^{\frac1t},
\end{align}
The lemma follows because if $Q \in \scriptQ_n$, then $\bigcup I_j^Q$ is covered by a bounded number of intervals in $\scriptD_n$, and each dyadic interval arises in only a bounded number of such coverings.  
\end{proof}

\begin{definition}
We define a family of Banach spaces $X^{p, q, r, s}$ with norms
$$
\|f\|_{X^{p, q, r, s}} := \bigl(\sum_n \bigl(\sum_{I \in \scriptD_n} 2^{n r (\frac1p - \frac1s)} \|f\|_{L^s (I)}^r \bigr)^{\frac qr} \bigr)^{\frac1q}.
$$
\end{definition}

Then Lemma~\ref{L:Xp,dt,dt,s} states that
$$
\|\scriptE f\|_q \lesssim \|f\|_{X^{p, d t, d t, s}},
$$
for $q > \tfrac{d^2 + d + 2}2$, $p = (\tfrac{2 q}{d^2 + d})'$, $t = (\tfrac qd)'$, $s = (\tfrac{2 q}{d^2})'$.  

\begin{lemma} \label{L:Xpqrs bound}
Assume that $1 < s < p < r \leq q < \infty$. Then $L^p \subseteq X^{p, q, r, s}$. Moreover, there exist $c_0 > 0$, $\theta > 0$ such that if $f \in L^p$ with $\|f\|_p = 1$, then
$$
\|f\|_{X^{p, q, r, s}} \lesssim \sup_{k \geq 0} \sup_{I} 2^{- c_0 k} \|f_I^k\|_p^{1 - \theta} \|f\|_p^\theta.
$$
Here the supremum is taken over all dyadic intervals $I$ and $f_I^k := f \rchi_I \rchi_{\{|f| \leq 2^k |I|^{- \frac1p}\}}$.  
\end{lemma}

We note that this immediately implies Proposition~\ref{P:refined Drury} in the range $p < d (\tfrac qd)'$, i.e.\ when $p < d + 2$.  

\begin{proof}[Proof of Lemma~\ref{L:Xpqrs bound}]
We will prove the superficially stronger estimate wherein we denote
$$
f^0_I := f \rchi_I \rchi_{\{|f| \leq |I|^{- \frac1p}\}}, \qquad f^k_I := f \rchi_I \rchi_{\{2^{k - 1} |I|^{- \frac1p} < |f| \leq 2^k |I|^{- \frac1p}\}}.
$$
Thus $f_I := f\rchi_I = \sum_{k \geq 0} f^k_I$. By H\"older's inequality,
\begin{align*}
\|f\|_{X^{p, q, r, s}}
& = \bigl( \sum_n \bigl(\sum_{I \in \scriptD_n} 2^{n r (\frac1p - \frac1s)} (\sum_{k \geq 0} \|f_I^k\|_s^s \bigr)^{\frac rs} \bigr)^{\frac qr} \bigr)^{\frac1q}
\\
& \leq \sup_n \sup_{I \in \scriptD_n} \sup_{k \geq 0} 2^{- c_0 k} (2^{n (\frac1p - \frac1s)} \|f_I^k\|_s \bigr)^{1 - \theta}
\\
& \qquad\qquad\qquad\times\bigl(\sum_n \bigl(\sum_{I \in \scriptD_n} 2^{n r (\frac1p - \frac1s) \theta} \bigl(\sum_{k \geq 0} 2^{c_0 k s}\|f_I^k\|_s^{s \theta} \bigr)^{\frac rs}\bigr)^{\frac qr}\bigr)^{\frac1q}.
\end{align*}
By H\"older's inequality, $2^{n (\frac1p - \frac1s)} \|f_I^k\|_s \leq \|f_I^k\|_p$, so it remains to bound the second term in the product on the right hand side.  

We begin with the $k = 0$ term. Since $r > p$, we may choose $\theta < 1$ sufficiently close to $1$ so that $r \theta > p > s$. Then using H\"older's inequality repeatedly and finally summing a geometric series,
\begin{align*}
&\bigl(\sum_n \bigl(\sum_{I \in \scriptD_n} 2^{n r \theta (\frac1p - \frac1s)} \|f_I^0\|_s^{\theta r} \bigr)^{\frac qr})^{\frac1q}
 \leq \bigl(\sum_n \bigl(\sum_{I \in \scriptD_n} 2^{n r \theta (\frac1p - \frac1{r \theta})} \|f_I^0\|_{r \theta}^{r \theta}\bigr)^{\frac qr} \bigr)^{\frac1q}
 \\
&\qquad \leq \bigl(\sum_n 2^{n r \theta (\frac1p - \frac1{r \theta})} \|f \rchi_{\{|f| \leq 2^{- \frac np}\}} \|_{r \theta}^{r \theta} \bigr)^{\frac1r}
\\
&\qquad = \bigl(\int |f|^{r \theta} \sum_{2^n < |f|^{- p}}2^{n r \theta (\frac1p - \frac1{r \theta})} \bigl)^{\frac1r} = \|f\|_p^{\frac pr}.
\end{align*}

Now we turn to the $k \geq 1$ terms. Let $c_0 < c_1 < c_2 < (p - s)\frac{\theta}{s}$. Then  several applications of H\"older's inequality and the triangle inequality give
\begin{align*}
&\bigl(\sum_n \bigl(\sum_{I \in \scriptD_n} 2^{n r (\frac1p - \frac1s) \theta}\bigl(\sum_{k \geq 1} 2^{c_0 k s} \|f_I^k\|_s^{s \theta} \bigr)^{\frac rs} \bigr)^{\frac qr}\bigr)^{\frac1q}
\\
&\qquad 
\lesssim \bigl(\sum_n \bigl(\sum_{I \in \scriptD_n} 2^{n r \theta (\frac1p - \frac1s)} \sum_{k \geq 1} 2^{c_1 k r}\|f_I^k\|_s^{r \theta}\bigr)^{\frac qr}\bigr)^{\frac1q}
\\
&\qquad 
\lesssim \bigl(\sum_{k \geq 1} 2^{c_2 k q} \sum_n 2^{n q \theta (\frac1p - \frac1s)} \bigl(\sum_{I \in \scriptD_n} \|f_I^k\|_s^{r \theta} \bigr)^{\frac qr}\bigr)^{\frac1q}
\\
&\qquad \leq \bigl(\sum_{k \geq 1} 2^{c_2 k q} \bigl(\sum_n 2^{n r \theta (\frac1p - \frac1s)} \sum_{I \in \scriptD_n} \|f_I^k\|_s^{r \theta} \bigr)^{\frac qr}\bigr)^{\frac1q}
\\
&\qquad \leq \bigl(\sum_{k \geq 1} 2^{c_2 k q} \bigl(\sum_n 2^{n s (\frac1p - \frac1s)} \sum_{I \in \scriptD_n} \|f_I^k\|_s^s \bigr)^{\frac{q \theta}s} \bigr)^{\frac1q}
\\
&\qquad \leq \bigl(\sum_{k \geq 1} 2^{c_2 k q} \bigl(\sum_n 2^{n s (\frac1p - \frac1s)} \int_{\{|f| \sim 2^{- \frac np} 2^k\}} |f|^s \bigr)^{\frac{q \theta} s}\bigr)^{\frac 1q}
\\
&\qquad \sim \bigl(\sum_{k \geq 1} 2^{c_2 k q} \bigl(2^{- k (p - s)} \|f\|_p^p \bigr)^{\frac{q \theta}s}\bigr)^{\frac1q}
 \sim \|f\|_p^{\frac{p \theta}{s}}.
\end{align*}
\end{proof}

In the case of larger $p$ (i.e., $d + 2 \leq p < \frac{d^2 + d + 2}{2}$), we will (roughly speaking) interpolate the bound in Proposition~\ref{P:refined Drury}, now established for sufficiently small $p$, with Drury's estimate $\|\scriptE f\|_q \lesssim \|f\|_p$. The details of this deduction are given in the next two lemmas. 

\begin{lemma} \label{L:SW}
Let $1 < p < \tfrac{d^2 + d + 2}2$ and $q = \tfrac{d (d + 1)}2 p'$. Let $f \in L^p$ and write $f = \sum_n 2^n f_n$, where $f_n := 2^{- n} f \rchi_{E_n}$, and $E_n := \{2^n \leq |f| < 2^{n + 1}\}$. Then   
$$
\|\scriptE f\|_q \lesssim \sup_n \|\scriptE 2^n f_n\|_q^\nu \|f\|_p^{1 - \nu},
$$
for some $0 < \nu  < 1$, depending only on $p$.  
\end{lemma}

\begin{proof}
We will prove the lemma by slightly adapting the proof of the Marcinkiewicz interpolation theorem from~\cite{SteinWeiss}. Write 
$$
\tfrac1p = \tfrac{1 - \theta}{p_0} + \tfrac\theta{p_1}, 
$$
for some $1 < p_0 < p < p_1 < \tfrac{d^2 + d + 2}2$; set $q_i := \tfrac{d (d + 1)}2 p_i'$, $i=0, 1$. Set 
$$
\nu := \min_{i = 0, 1} \frac{1/p_i - 1/q_i}{1 + 1/p_i - 1/q_i}.
$$
Note that $0 < \nu < 1$.  

Let $g \in L^{q'}$ and decompose $g$ analogously to $f$: $g = \sum_n 2^n g_n$, where $g_n = 2^{- n} g \rchi_{F_n}$, $F_n := \{2^n \leq |g| < 2^{n + 1}\}$. We may assume that $\|f\|_p = \|g\|_{q'} = 1$ and thus $\sum_n 2^{n p} |E_n| \sim \sum_m 2^{m q'} |F_m| \sim 1$.  

By H\"older's inequality, Drury's theorem, and the definition of the $f_n, g_m$,
\begin{align*}
& \langle \scriptE f, g\rangle 
\lesssim  \sum_{n, m} \langle 2^n \scriptE f_n, 2^m g_m\rangle
\\
&\qquad \leq \sum_{n, m} \|2^n \scriptE f_n\|_q^\nu \|2^m g_m\|_{q'}^\nu \min_{i = 0, 1} \|2^n \scriptE f_n\|_{q_i}^{1 - \nu} \|2^m g_m\|_{q_i'}^{1 - \nu}
\\
&\qquad \lesssim \sup_n \|2^n \scriptE f_n\|_q^\nu \sum_{n, m} \min_{i = 0, 1} \|2^n f_n\|_{p_i}^{1 - \nu} \|2^m g_m\|_{q_i'}^{1 - \nu}
\\
&\qquad \lesssim \sup_n \|2^n \scriptE f_n\|_q^\nu \sum_{n, m} (2^{n + m} \min_{i =0, 1} |E_n|^{\frac1{p_i}} |F_m|^{\frac 1{q_i'}})^{1 - \nu}.
\end{align*}
It remains to bound the sum on the right side of this inequality.

Were it the case that $2^n |E_n|^{\frac1p} = 2^m |F_m|^{\frac1{q'}} = 1$, 
$$
|E_n|^{\frac1{p_0}} |F_m|^{\frac1{q_0'}} \leq |E_n|^{\frac1{p_1}} |F_m|^{\frac1{q_1'}}
$$
would hold if and only if $n A \leq - m B$, where
$$
A := p (\tfrac1{p_1} - \tfrac1{p_0}), \qquad B := q' (\tfrac1{q_1'} - \tfrac1{q_0'}).
$$
In any case, 
\begin{align}\notag
&\sum_{n, m} (2^{n + m} \min_{i = 0, 1} |E_n|^{\frac1{p_i}} |F_m|^{\frac 1{q_i'}})^{1 - \nu}
\\ \label{E:nA+mB}
&\quad \leq \sum_{n A + m B \leq 0}(2^{n + m} |E_n|^{\frac1{p_0}} |F_m|^{\frac 1{q_0'}})^{1 - \nu} + \sum_{n A + m B > 0}(2^{n + m} |E_n|^{\frac1{p_1}} |F_m|^{\frac 1{q_1'}})^{1 - \nu}
\end{align}

We begin with the first summand on the right of \eqref{E:nA+mB}. Simple arithmetic, followed by H\"older's inequality (since $\frac{1 - \nu}{p_0} + \frac{1 - \nu}{q_0'} \geq 1$) gives
\begin{align*}
&\sum_{n A + m B \leq 0}(2^{n + m} |E_n|^{\frac1{p_0}} |F_m|^{\frac 1{q_0'}})^{1 - \nu}
\\
&\qquad = \sum_{n A + m B \leq 0} 2^{\theta (n A + m B)(1 - \nu)}(2^{n p} |E_n|)^{\frac{1 - \nu}{p_0}}(2^{m q'} |F_m|)^{\frac{1 - \nu}{q_0'}}
\\
&\qquad \leq \sum_{k \leq 0} 2^{\theta k (1 - \nu)} \sum_n (2^{n p} |E_n|)^{\frac{1 - \nu}{p_0}} \sum_{\ceil{n A + m B} = k}(2^{m q'} |F_m|)^{\frac{1 - \nu}{q_0'}}
\\
&\qquad = \sum_{k \leq 0} 2^{\theta k (1 - \nu)} \sum_n (2^{n p} |E_n|)^{\frac{1 - \nu}{p_0}} \sum_{\frac{k - n A - 1}{B} < m \leq \frac{k - n A}{B}}(2^{m q'} |F_m|)^{\frac{1 - \nu}{q_0'}}
\\
&\qquad \lesssim \sum_{k \leq 0} 2^{\theta k (1 - \nu)} \big(\sum_n 2^{n p} |E_n| \big)^{\frac{1 - \nu}{p_0}}  \big(\sum_n \sum_{\frac{k - n A - 1}{B} < m \leq \frac{k - n A}{B}}2^{m q'} |F_m| \big)^{\frac{1 - \nu}{q_0'}}\lesssim 1.
\end{align*}

The proof of Lemma~\ref{L:SW} is complete modulo the bound for the second term on the right of \eqref{E:nA+mB}, which can be proved in an analogous fashion.  
\end{proof}

\begin{lemma} \label{L:refined char fxns}
Proposition~\ref{P:refined Drury} holds for functions $|f| \sim \lambda \rchi_E$, for $\lambda > 0$ and $E$ a measurable subset of $\R$, with bounds independent of $\lambda, E$.    
\end{lemma}

\begin{proof}
We may assume that $\lambda = 1$. Choose $p_0, p_1, q_0, q_1, \theta$ as in the proof of Lemma~\ref{L:SW}, with the additional assumption that $p_0 < d+2$. By H\"older's inequality, then the remark following Lemma~\ref{L:Xpqrs bound} and Drury's theorem,
$$
\|\scriptE f\|_q \leq \|\scriptE f\|_{q_0}^{1 - \theta} \|\scriptE f\|_{q_1}^\theta \lesssim (\sup_{k \in \Z} \sup_{I \in \scriptD_k} \sup_{n \geq 0} 2^{- c_{p_0} n} \|(\rchi_E)_I^n\|_{p_0})^{\theta'} \|\rchi_E\|_{p_0}^{1 - \theta - \theta'} \|\rchi_E\|_{p_1}^\theta,
$$
where $0 < \theta' < 1 - \theta$ arises by applying Proposition~\ref{P:refined Drury} with exponents $p_0, q_0$ (where it has already been established).  

The proof is now just a matter of unwinding the definition of $(\rchi_E)_I^n$ and performing some arithmetic. Observe that
$$
\|(\rchi_E)_I^n\|_{p_0} = \begin{cases} |E \cap I|^{\frac1{p_0}}, \: \text{if} \: 1 < 2^n |E|^{\frac1{p_0}} |I|^{- \frac1{p_0}}, \\ 0, \: \text{otherwise},\end{cases}
$$
and analogously with $p$ in place of $p_0$.  Thus
we may rewrite
\begin{align*}
\sup_{k \in \Z} \sup_{I \in \scriptD_k} \sup_{n \geq 0} 2^{- c_{p_0} n} \|(\rchi_E)_I^n\|_{p_0} 
& = \sup_{k \in \Z} \sup_{I \in \scriptD_k} \min\{1, \bigl(\tfrac{|E|}{2^k} \bigr)^{\frac{c_{p_0}}{p_0}}\} |E \cap I|^{\frac1{p_0}}
\\
& \sim \bigl(\sup_{k \in \Z} \sup_{I \in \scriptD_k} \sup_{n \geq 0} 2^{- c_{p_0} n} \|(\rchi_E)_I^n\|_{p} \bigr)^{\frac p{p_0}}.
\end{align*}
Finally, we obtain
\begin{align*}
\|\scriptE f\|_q &\lesssim (\sup_{k \in \Z} \sup_{I \in \scriptD_k} \sup_{n \geq 0} 2^{- c_{p_0} n} \|(\rchi_E)_I^n\|_p)^\vartheta \|\rchi_E\|_p^{1 - \vartheta}
\\
& \sim (\sup_{k \in \Z} \sup_{I \in \scriptD_k} \sup_{n \geq 0} 2^{- c_{p_0} n} \|f_I^n\|_p)^\vartheta \|f\|_p^{1 - \vartheta},
\end{align*}
where $\vartheta = \tfrac{\theta' p}{p_0}$. Note that $0 < \vartheta < \tfrac{(1 - \theta) p}{p_0} < 1$.  
\end{proof}

Proposition~\ref{P:refined Drury} in the cases $p \geq d+2$ follows by first applying Lemma~\ref{L:SW}, then applying Lemma~\ref{L:refined char fxns} to the supremum term in the conclusion of Lemma~\ref{L:SW}, and finally observing that, in the decomposition in Lemma~\ref{L:SW}, 
$$
|(2^m f_m)_I^n| \leq |(f)_I^n|,
$$
for all integers $m, n$ and intervals $I$.  
\end{proof}

The following proposition is somewhat easier to use, though it only applies in a more limited range of exponents. The proof requires only a small modification in the argument leading to Proposition~\ref{P:refined Drury}.  

\begin{proposition} \label{P:refined Drury 2}
Let $1 < p < d + 2$, and let $q := \tfrac{d (d + 1)}2 p'$. There exists $0 < \theta = \theta_p < 1$ such that for $f \in L^p$,
\begin{equation} \label{E:refined Drury 2}
\|\scriptE f\|_q \lesssim \bigl(\sup_I |I|^{- \frac1{p'}} \|\scriptE f_I\|_\infty \bigr)^{1 - \theta} \|f\|_p^\theta.
\end{equation}
Here, the supremum is taken over dyadic intervals $I$, and $f_I := f \rchi_I$.
\end{proposition}
\begin{proof}
Let $q_1 < q$ with $\theta := \tfrac{q_1}q$ sufficiently close to $1$ for later purposes. By \eqref{E:Xp,dt,dt,s}, H\"older, Lemma~\ref{L:scaled multilinear extn}, arithmetic, and another application of H\"older,
\begin{align*}
\|\scriptE f\|_q & \lesssim \bigl(\sum_n \sum_{Q = \prod I_j \in \scriptQ_n} \|\prod_{j = 1}^d \scriptE f_{I_j}\|_{\frac qd}^t \bigr)^{\frac1{t d}}
\\
& \lesssim \bigl(\sum_n \sum_{Q = \prod I_j \in \scriptQ_n} \|\prod_{j = 1}^d \scriptE f_{I_j}\|_\infty^{t (1 - \theta)} \|\prod_{j = 1}^d \scriptE f_{I_j}\|_{\frac {q_1}d}^{t \theta} \bigr)^{\frac1{t d}}
\\ 
& \lesssim \bigl(\sum_n \sum_{Q = \prod I_j \in \scriptQ_n} \|\prod_{j = 1}^d \scriptE f_{I_j}\|_\infty^{t (1 - \theta)} 2^{n d t \theta(\frac1{p_1} - \frac1{s_1})} \prod_{j = 1}^d \|f_{I_j}\|_{s_1}^{t \theta}\bigr)^{\frac1{t d}}
\\  
& \leq \bigl(\sum_n \sum_{Q = \prod I_j \in \scriptQ_n} (2^{n d (\frac1p - 1)} \max_{j = 1, \ldots, d} \|\scriptE f_{I_j}\|_\infty^d)^{t (1 - \theta)} 2^{n d t \theta (\frac1p - \frac1{s_1})} \prod_{j = 1}^d \|f_{I_j}\|_{s_1}^{t \theta} \bigr)^{\frac1{t d}}
\\ 
& \leq \bigl(\sup_I |I|^{\frac1p-1} \|\scriptE f_I\|_\infty\bigr)^{1 - \theta} \|f\|_{X^{p, d t \theta, d t \theta, s_1}}^\theta.
\end{align*}
Here $s_1 := (\frac{2 q_1}{d^2})', t := (\frac{q}{d})'$. For $p < d + 2$ and $\theta$ sufficiently close to $1$, $1 < s_1 < p < d t \theta < \infty$, so an application of Lemma~\ref{L:Xpqrs bound} completes the proof.
\end{proof}

\section{Frequency Localization} \label{S:Freq loc}
In this section we prove that any near extremizer of $\scriptE$ is uniformly bounded and is supported on a compact set around the origin possibly after applying symmetry, if we allow ourselves to lose a small amount of $L^p$-mass. Below is the precise statement.

\begin{proposition}\label{P:freq loc}
Let $1 < p < \tfrac{d^2 + d + 2}2$, and let $q := \tfrac{d (d + 1)}2 p'$. For each $\epsilon > 0$, there exist $\delta > 0$ and $R < \infty$ such that for each nonzero function $f$ satisfying $\|\scriptE f\|_q \geq B_p (1 - \delta) \|f\|_p$, there exists a symmetry $S$ such that the following holds.  
$$
\|Sf\|_{L^p \big(\{|t| > R\} \cup \{|Sf| > R \|f\|_p\}\big)} < \epsilon \|f\|_{p}.
$$
\end{proposition}

We start with the following lemma.
\begin{lemma} \label{L:nibble}
Let $1 < p < \tfrac{d^2 + d + 2}2$, and $q := \tfrac{d (d + 1)}2 p'$. There exists a sequence $\rho_k \to 0$ such that for every $f \in L^p(\R)$, there exists a sequence $\{I_k\}$ of dyadic intervals such that if $\{f^{> k}\}$ is inductively defined by
\begin{equation}\label{partition:f}
f^{> 0} := f, \qquad f^k := f^{> k - 1} \rchi_{\{|f| < 2^k |I_k|^{- \frac1p} \|f\|_p\}} \rchi_{I_k}, \qquad f^{> k} := f^{> k - 1} - f^k,
\end{equation}
then for any measurable function $h^{>k}$ with $|h^{> k}| = \rchi_E |f^{> k}|$, for some measurable set $E$,
$$
\|\scriptE h^{> k}\|_q \leq \rho_k \|f\|_p.  
$$
\end{lemma}

\begin{proof}
Let $0 \not \equiv f  \in L^p$. Multiplying by a constant if needed, we may assume that $\|f\|_p = 1$. By the Dominated Convergence Theorem, given $f^{> k - 1}$, we choose a dyadic $I_k$ to maximize $\|f^k\|_p$. With the sequence $\{f^k\}$ and $\{f^{> k}\}$ as defined in~\ref{partition:f}, let $K \in \N$ and set
$$
A_K := \sup_{I\,\text{dyadic}} \|(f^{> 2 K})_I^K\|_p.
$$
By the maximality property of the $I_k$, $A_K \leq \|f^{K + j}\|_p$, for each $0 \leq j \leq K$. Hence by the disjointness of the supports of the $f^k$'s,
$$
KA_K^p \leq \sum_{j = 1}^K \|f^{K + j}\|_p^p \leq 1.
$$
By Proposition~\ref{P:refined Drury}, for any measurable function $|h^{> 2 K}| =\rchi_E |f^{> 2 K}|$,
$$
\|\scriptE h^{> 2 K}\|_q \lesssim \max\{2^{- c_p \theta_p K}, \tfrac1{K^{\frac{\theta}{p}}}\} \lesssim K^{- \frac{\theta_p}{p}}.
$$
This completes the proof with $\rho_k := C (\frac k2)^{- \frac{\theta_p}{p}}$ with $\theta_p$ as in Proposition~\ref{P:refined Drury} and $C$ a sufficiently large constant.
\end{proof}

\begin{proof}[Proof of Proposition~\ref{P:freq loc}]
Since we can always take $\delta < \tfrac12$, it suffices to consider those $f$ for which $\|f\|_p = 1$ and $\|\scriptE f\|_q > \tfrac12 B_p \|f\|_p$. Let $I_{k, f}$ denote the dyadic intervals from Lemma~\ref{L:nibble} and set 
$$
f^{\leq k} := f - f^{> k}, \qquad\text{and so}\,\, f^k = f^{> k - 1} - f^{> k}.
$$
By Lemma~\ref{L:nibble}, $\|f^{k_f}\|_p \gtrsim 1$ for some $k_f \lesssim 1$. Applying a symmetry if needed, we may assume that $I_{k_f, f}$ is the unit interval. We will prove (under these assumptions on $f$) that the conclusion of the proposition holds with $S$ equal to the identity.  

If the conclusion were to fail, there would exist some $\eps > 0$ and sequence of functions $\{f_n\} \subseteq L^p$ with $\|f_n\|_p = 1$, $k_n := k_{f_n} \lesssim 1$, $\|f_n^{k_n}\|_p \gtrsim 1$, $I_{k_n, f_n} = [0, 1]$, $\|\scriptE f_n\|_q \geq B_p (1 - n^{- 1})$, and 
$$
\|f_n\|_{L^p \big(\{|t| > n\} \cup \{|f_n| > n\}\big)} > \eps.
$$

By pigeonholing and passing to a subsequence, we may assume that $k_n = k_0 \lesssim 1$ for all $n$. Write 
$$
I_n^k := I_{k, f_n} =: \xi_n^k + [0, \ell_n^k], \qquad k, n \in \N.
$$
Passing to a subsequence, we may assume that $\{\xi_n^k\}$ and $\{\ell_n^k\}$ converge in $[- \infty, \infty]$ and $[0, \infty]$, respectively, for each $k$ and that $\{\|f_n^k\|_p\}$ converges for all $k$. For each $k$, say that the $k$ is negligible if $\|f_n^k\|_p \to 0$, that the $k$ is good if it is negligible or if $\{\xi_n^k\}$ converges in $\R$ and $\{\ell_n^k\}$ converges in $(0,\infty)$. Say that the $k$ is bad if it is not good. For bad $k$: say $k$ is \textit{long} if $\ell_n^k \to \infty$ as $n \to \infty$, \textit{short} if $\ell_n^k \to 0$ as $n \to \infty$; otherwise it must be \textit{far}, i.e.\ $|\xi_n^k| \to \infty$.

We will prove that every $k$ is good. Assuming this for now, we complete the proof of the proposition. Since every $k$ is good, for each fixed $K$ and sufficiently large $n$ (depending on $K$),
$$
|f_n| \rchi_{(\{|t| > n\} \cup \{|f_n| > n\})} \leq |f_n^{> K}|.
$$
Therefore
$$
\liminf_{n \to \infty} \|f_n^{> K}\|_p > \eps,
$$ 
for every $K$.  
Since the supports of the $f_n^{\leq K}$ and $f_n^{> K}$ intersect on a set of measure zero, 
$$
\limsup_{n \to \infty} \|f_n^{\leq K}\|_p = \bigl(\|f_n\|_p^p - \|f_n^{> K}\|_p^p)^{\frac1p} < (1 - \eps^p)^{\frac1p}.
$$
We conclude that $\|\scriptE f_n^{\leq K}\|_q < B_p (1 - \eps^p)^{\frac1p}$. On the other hand, by Lemma~\ref{L:nibble}, $\|\scriptE f_n^{ > K}\|_q \leq \rho_K$. Therefore $\|\scriptE f_n^{\leq K}\|_q \geq B_p - \rho_K$. For $K$ sufficiently large, depending on $\eps$, 
$$
B_p (1 - \eps^p)^{\frac1p} > B_p - \rho_K,
$$
leading to a contradiction.  

It remains to prove that every $k$ is good. Define 
$$
g_n^{\leq K} := \sum_{k \leq K, \text{\,good}} f_n^k, \qquad b_n^{\leq K} := \sum_{k \leq K, \text{\,bad}} f_n^k.  
$$
Then $f_n^{\leq K} = g_n^{\leq K} + b_n^{\leq K}$, and the supports of $g_n^{\leq K}$ and $b_n^{\leq K}$ have measure zero intersection and thus obey the $L^p$ orthogonality condition
$$
\|f_n^{\leq K}\|_p^p = \|g_n^{\leq K}\|_p^p + \|b_n^{\leq K}\|_p^p.
$$
Our next lemma shows that we have $L^q$-orthogonality in limit.  

\begin{lemma} \label{L:Lq freq orthog}
For every $1 \leq K < \infty$,
$$
\lim_{n \to \infty} \|\scriptE f_n^{\leq K}\|_q^q -\bigl(\|\scriptE g_n^{\leq K}\|_q^q + \|\scriptE b_n^{\leq K}\|_q^q\bigr) = 0.
$$
\end{lemma}

Assuming Lemma~\ref{L:Lq freq orthog}, we complete the proof of the proposition by showing that every $k$ is good. Suppose, by way of contradiction, that $k_0$ is bad. By definition, $\limsup_{n \to \infty} \|f_n^{k_0}\|_p > \eps$, for some $\eps > 0$. Passing to a subsequence, we may assume that $\|f_n^{k_0}\|_p > \eps$ for all $n$. Therefore $\|b_n^{\leq K}\|_p > \eps$ for all $K \geq k_0$ and all $n$. Taking a smaller $\eps$ if needed, $\|g_n^{\leq K}\|_p \geq \|f_n^{k_0}\|_p > \eps$, for all $K \geq k_0$ and all $n$. We may further assume that $\eps^p < \tfrac12$. We can use these $L^p$ estimates to bound the extension for sufficiently large $K$:
\begin{align*}
&\limsup_{n \to \infty} \|\scriptE f_n^{\leq K}\|_q^q 
\leq \limsup_{n \to \infty} \|\scriptE g_n^{\leq K}\|_q^q + \|\scriptE b_n^{\leq K}\|_q^q
\\
&\qquad \leq \limsup_{n \to \infty} B_p^q \bigl(\|g_n^{\leq K}\|_p^q + \|b_n^{\leq K}\|_p^q\bigr) \leq \bigl[(1 - \eps^p)^{\frac qp} + \eps^q\bigr] B_p^q.
\end{align*}
For the last inequality, we have used that for $g, b \geq 0$, $g^p + b^p \leq 1$, $b, g > \eps$, and $\eps^{p} < \frac{1}{2}$,
$$
g^{q}+b^{q} \leq (1 - \eps^p)^{\frac qp} + \eps^q,
$$
which in turn follows from basic calculus.  

Crucially, $\bigl[(1 - \eps^p)^{\frac qp} + \eps^q\bigr] =: c_\eps < 1$. On the other hand, by our hypothesis on $\{f_n\}$ and Lemma~\ref{L:nibble}, 
$$
B_p = \lim_{n \to \infty} \|\scriptE f_n\|_q = \lim_{K\to\infty}\lim_{n \to \infty} \|\scriptE f_n^{\leq K}\|_q.  
$$
\end{proof}

It remains to prove Lemma~\ref{L:Lq freq orthog}.
\begin{proof}[Proof of Lemma~\ref{L:Lq freq orthog}]
By elementary calculus, for all $q > 1$ and $a, b \geq 0$,
$$
|(a + b)^q - a^q - b^q| \lesssim_{q} \bigl(a b^{q - 1} + b a^{q - 1}\bigr).
$$
Since $q > \frac{d^2 + d + 2}{2} > 2$, we may apply this inequality to our good and bad part of the function to see that  
\begin{align*}
&\|\scriptE f_n^{\leq K}\|_q^q - \|\scriptE g_n^{\leq K}\|_q^q - \|\scriptE b_n^{\leq K}\|_q^q
\\
&\qquad \lesssim_{q, K} \sum_{k \leq K \, \text{good}} \sum_{k' \leq K\, \text{bad}} \int |\scriptE f_n^k| |\scriptE f_n^{k'}|^{q - 1} + |\scriptE f_n^k|^{q - 1} |\scriptE f_n^{k'}|.
\end{align*}
Moreover, by H\"older's inequality, and boundedness of the $\scriptE f_n^j$,
\begin{align*}
\int |\scriptE f_n^k| |\scriptE f_n^{k'}|^{q - 1} + |\scriptE f_n^k|^{q - 1} |\scriptE f_n^{k'}| & \leq \bigl(\|\scriptE f_n^k\|_q^{q - 2} + \|\scriptE f_n^{k'}\|_q^{q - 2}\bigr) \|\scriptE f_n^k \scriptE f_n^{k'}\|_{\frac{q}{2}}
\\
&\lesssim \|\scriptE f_n^k \scriptE f_n^{k'}\|_{\frac{q}{2}},
\end{align*}
so it suffices to prove that $\|\scriptE f_n^k \scriptE f_n^{k'}\|_{\frac{q}{2}} \to 0$ when $k$ is good and $k'$ is bad.  

To this end, choose $q_+ < q < q_-$ with $q_{\pm} = \tfrac{d (d + 1)}2 p_{\pm}'$ and $\frac{2}{q} = \frac{1}{q_+} + \frac{1}{q_-}$. Therefore $p_- < p < p_+$. Because $k$ is good, $f_n^k$ remains bounded in $L^{p_\pm}$ (in fact, in every Lebesgue space) as $n \to \infty$. When the $k'$ is short, $\|f_n^{k'}\|_{p_-} \rightarrow 0$. Therefore
$$
\|\scriptE f_n^k \scriptE f_n^{k'}\|_{\frac{q}{2}} \leq \|\scriptE f_n^k\|_{q_+} \|\scriptE f_n^{k'}\|_{q_-} \lesssim \|f_n^{k'}\|_{p_-} \rightarrow 0.
$$
Similarly when the $k'$ is long, $\|f_n^{k'}\|_{p_+} \rightarrow 0$ and thus
$$
\|\scriptE f_n^k \scriptE f_n^{k'}\|_{\frac{q}{2}} \lesssim \|f_n^{k'}\|_{p_+} \rightarrow 0.
$$

Finally, suppose the $k'$ is far (and neither short nor long). We set $L := \lim l_n^k + \lim l_n^{k'}$, and assume that $n$ is sufficiently large so that $l_n^k + l_n^{k'} \leq 2 L$ and $|\xi_n^k - \xi_n^{k'}| \geq 100 L$. By the arithmetic-geometric mean inequality, Lemma~\ref{L:scaled multilinear extn}, and H\"older's inequality with $s := (\tfrac{2 q}{d^2})' < p$,
\begin{align*}
&\|\scriptE f_n^k \scriptE f_n^{k'}\|_{\frac q2} \leq \|(\scriptE f_n^k)^{d - 1} \scriptE f_n^{k'}\|_{\frac qd} + \|\scriptE f_n^k (\scriptE f_n^{k'})^{d - 1}\|_{\frac qd}
\\
&\qquad \lesssim |\xi_n^k - \xi_n^{k'}|^{- d (\frac1s - \frac1p)} (\|f_n^k\|_s^{d - 1}\|f_n^{k'}\|_s + \|f_n^k\|_s \|f_n^{k'}\|_s^{d - 1})
\\
&\qquad \lesssim \bigl(\tfrac L{|\xi_n^k - \xi_n^{k'}|})^{d (\frac1s - \frac1p)} \to 0.
\end{align*}
\end{proof}

\section{A profile decomposition for frequency localized sequences}

\begin{proposition}\label{P:profile}
Let $q = \tfrac{d (d + 1)}2 p' > p > 1$. Let $\{f_n\}$ be a sequence of measurable functions with $\supp f_n \subseteq [- R, R]$ and $|f_n| \leq R$ for all $n$. Then, after passing to a subsequence, there exist $\{x_n^j\}_{n, j \geq 1} \subseteq \R^d$ and $\{\phi^j\} \subseteq L^p$ such that the following hold with
$$
w_n^J := f_n - \sum_{j = 1}^J e^{- i x_n^j \cdot \gamma} \phi^j.
$$
\begin{enumerate}[\rm(i)]
\item $\lim_{n \to \infty} |x_n^j - x_n^{j'}| = \infty$, for all $j \neq j'$;
\item $e^{i x_n^j \cdot \gamma} f_n \rightharpoonup \phi^j$, weakly in $L^p$, for all $j$;
\item $\lim_{n \to \infty} \|\scriptE f_n\|_q^q - \sum_{j = 1}^J \|\scriptE \phi^j\|_q^q - \|\scriptE w_n^J\|_q^q = 0$ for all $J$;
\item $\lim_{J \to {\infty}} \lim_{n \to \infty} \|\scriptE w_n^J\|_q = 0$;
\item $\bigl(\sum_{j = 1}^{\infty} \|\phi^j\|_p^{{\tilde{p}}}\bigr)^{\frac1{\tilde p}} \leq \liminf_{n \to \infty} \|f_n\|_p$, where $\tilde p := \max\{p, p'\}$.  
\end{enumerate}
\end{proposition}

\begin{remark}
Here we allow for the possibility that the $\phi^j$ with $j$ sufficiently large might all be identically zero.   
\end{remark}

The essential step is finding a nonzero weak limit. 

\begin{lemma}\label{L:nonzero wk lim}
Let $q = \tfrac{d (d + 1)} 2 p' > p > 1$ with $p < d + 2$. There exists $C > 0$ such that for any sequence $\{f_n\} \subseteq L^p$ with $|f_n| \leq R \rchi_{[- R, R]}$,
$$
\|f_n\|_p \leq A, \qtq{and} \|\scriptE f_n\|_q \geq \eps,
$$
there exists a sequence $\{x_n\} \subseteq \R^d$ such that, after passing to a subsequence, $e^{ix_n \cdot \gamma} f_n \rightharpoonup \phi$, weakly in $L^p$, for some $\phi \in L^p$ with $\|\phi\|_p \gtrsim \eps (\tfrac \eps A)^C$. Here the implicit constant is independent of $R$.  
\end{lemma}

\begin{proof}
By \eqref{E:refined Drury 2}, for each $n$, there exists a dyadic interval $I_n$ such that
\begin{equation} \label{E:first eps A lb}
\eps \leq \bigl(|I_n|^{- \frac1{p'}} \|\scriptE(f_n)_{I_n}\|_\infty\bigr)^\theta A^{1 - \theta}.
\end{equation}
By H\"older's inequality, 
$$
|I_n|^{- \frac1{p'}} \|\scriptE(f_n)_{I_n}\|_\infty \leq C_R \min\{|I_n|^{- \frac1{p'}}, |I_n \cap [- R, R]|^{\frac1p}\}.
$$
Hence, in the terminology of the previous section, $\{I_n\}$ cannot be long, short, nor far, and so, after passing to a subsequence, we may assume that $I_n = I$ is independent of $n$.   

By \eqref{E:first eps A lb}, for each $n$, there exists $x_n \in \R^d$ such that
$$
\eps \bigl(\tfrac\eps A)^{\frac{1 - \theta}{\theta}} \lesssim |I|^{- \frac1{p'}} |\scriptE(f_n)_I(x_n)| = |I|^{- \frac1{p'}} |\scriptE(e^{i x_n \cdot \gamma} f_n)_I(0)|.
$$
By boundedness of the sequence $\{f_n\}$, after passing to a subsequence, $e^{i x_n \gamma} f_n \rightharpoonup \phi$, weakly in $L^p$, for some $\phi \in L^p$. Along this same subsequence, we then have $e^{i x_n \gamma} (f_n)_I \rightharpoonup \phi_I$. By compactness of their support, 
$$
\scriptE \phi_I(0) = \lim_{n \to \infty} \scriptE(e^{i x_n \cdot \gamma} f_n)_I(0).
$$
Therefore by H\"older's inequality,
$$
\eps \bigl(\tfrac\eps A)^{\frac{1 - \theta}{\theta}} \lesssim |I|^{- \frac1{p'}} |\scriptE \phi_I(0)| \leq |I|^{- \frac1{p'}} \|\phi_I\|_1 \leq \|\phi_I\|_p.
$$
\end{proof}

The remainder of the section is devoted to the proof of Proposition~\ref{P:profile}. We prove the proposition first in the case $p = 2$, and then in the general case.

\begin{proof}[Proof of Proposition~\ref{P:profile} when $p = 2$]
In the case $p = 2, q = q_2 := d (d + 1)$, we may replace (v) with the stronger condition that for all $J$,
$$
\qtq{(v')} \lim_{n \to \infty} \|f_n\|_2^2 - \sum_{j = 1}^J \|\phi^j\|_2^2 - \|w_n^J\|_2^2 = 0.
$$

Suppose that we are given $1 \leq J_1 < \infty$ and sequences $\{x_n^j\}_{n \in \N, j < J_1} \subseteq \R^d$, $\{\phi^j\}_{j < J_1} \subseteq L^2$ such that (i - iii) and (v') hold for all $J < J_1$. If $\lim \|\scriptE w_n^{J_1 - 1}\|_{q_2} = 0$, we are done after setting $\phi^j \equiv 0$ for $j \geq J_1$. Otherwise, after passing to a subsequence, for sufficiently large $n$, $\|\scriptE w_n^{J_1 - 1}\|_{q_2} > \eps > 0$. By (ii), $|w_n^{J_1 - 1}| \leq J_1 R \rchi_{[- R, R]}$, and by (v'), 
$$
\limsup \|w_n^{J_1 - 1}\|_2 \leq \limsup \|f_n\|_2 =: A.  
$$  
Therefore, by Lemma~\ref{L:nonzero wk lim}, there exists $\{x_n^{J_1}\} \subseteq \Rd$ and a subsequence along which
$$
e^{i x_n^{J_1} \cdot \gamma} w_n^{J_1 - 1} \rightharpoonup \phi^{J_1}\qquad \text{weakly in}\, L^2,
$$
with $\|\phi^{J_1}\|_2 \gtrsim \eps^C$. This immediately implies that
$$
\lim_{n \to \infty} \|e^{i x_n^{J_1} \cdot \gamma} w_n^{J_1 - 1}\|_2^2 - \|\phi^{J_1}\|_2^2 - \|e^{i x_n^{J_1} \cdot \gamma} w_n^{J_1 - 1} - \phi^{J_1}\|_2^2 = 0,
$$
so (v') holds with $J = J_1$. By the compact support condition, $\scriptE( e^{i x_n^{J_1} \cdot \gamma} w_n^{J_1 - 1}) \to \scriptE \phi^{J_1}$, a.e., so by the Brezis--Lieb lemma,
$$
\lim_{n \to \infty} \|\scriptE (e^{i x_n^{J_1} \cdot \gamma} w_n^{J_1 - 1})\|_q^q - \|\scriptE \phi^{J_1}\|_q^q - \|\scriptE (e^{i x_n^{J_1} \cdot \gamma} w_n^{J_1 - 1} - \phi^{J_1})\|_q^q = 0.
$$
Therefore (iii) holds with $J = J_1$.  

Suppose that $|x_n^{J_1} - x_n^{j_0}| \not \to \infty$, for some $j_0 < J_1$. Passing to a subsequence, we may assume that $x_n^{J_1} - x_n^{j_0} \to x$. Then multiplication by $e^{i (x_n^{J_1} - x_n^{j_0}) \cdot \gamma}$ converges to multiplication by $e^{i x \cdot \gamma}$ in the strong operator topology. Moreover, by (i), multiplication by $e^{i (x_n^{J_1} - x_n^{j}) \cdot \gamma}$, $j \neq j_0$, converges to zero in the weak operator topology. Therefore
\begin{align*}
\wklim e^{i x_n^{J_1} \cdot \gamma} w_n^{J_1 - 1} 
& = \wklim e^{i (x_n^{J_1} - x_n^{j_0}) \cdot \gamma} (e^{i x_n^{j_0}  \cdot \gamma} f_n - \phi^{j_0}) + \sum_{j \neq j_0} e^{i (x_n^{J_1} - x_n^j) \cdot \gamma} \phi^j = 0.
\end{align*}
On the other hand, the left hand side of the preceding equals $\phi^{J_1}$, which is nonzero, a contradiction. Tracing back, (i) holds for indices bounded by $J_1$. The proof of (ii) is similar: Along our subsequence,
 $$
 \wklim e^{i x_n^{J_1} \cdot \gamma} f_n = \wklim \sum_{j = 1}^{J_1 - 1} e^{i (x_n^{J_1} - x_n^j) \cdot \gamma} \phi^j + e^{i x_n^{J_1} \cdot \gamma} w_n^{J_1 - 1} = \sum_{j = 1}^{J_1 - 1} 0 + \phi^{J_1}.
 $$
 
It remains to verify (iv). Let $\eps^J := \limsup_{n \to \infty} \|\scriptE w_n^J\|_q$. If $\eps^J > \eps > 0$ infinitely often, then, as we have seen, $\|\phi^J\|_2 \gtrsim \eps^C$ infinitely often, which is impossible by (v').  
\end{proof}

Having proved Proposition~\ref{P:profile} in the case $p = p_2 := 2$, $q = q_2 := d (d + 1)$, we turn to the general case.  

\begin{proof}[Proof of Proposition~\ref{P:profile} for $d^2 + d \neq q = \tfrac{d^2 + d}2 p' > p > 1$]
Fix an exponent $\infty > q_1 >  \tfrac{d^2 + d + 2}2$ so that $q$ lies strictly between $q_1$ and $q_2$. As $\{f_n\}$ is bounded in $L^2$, we may apply the $L^2$-based profile decomposition to determine $\{x_n^j\}, \{\phi^j\}$. That (i) and (ii) hold are immediate. That (iii) holds follows from the Brezis--Lieb lemma and the above argument in the case $q = q_2$. The Brezis--Lieb lemma also implies that (iii) holds in the case $q = q_1$, so $\|\scriptE w_n^J\|_{q_1}$ is bounded, uniformly in $J$ (albeit with a constant that depends on $R$). Choosing $\theta$ so that $\tfrac1q = \tfrac\theta{q_1} + \tfrac{1 - \theta}{q_2}$,
$$
\lim_{J \to \infty} \lim_{n \to \infty} \|\scriptE w_n^J\|_q \lesssim_R \lim_{J \to \infty} \lim_{n \to \infty} \|\scriptE w_n^J\|_{q_2}^{1 - \theta} = 0,
$$
so (iv) holds as well.

It remains to prove (v). Let us fix $J$ and let $\epsilon > 0$. We choose compactly supported smooth nonnegative functions $a, b$ satisfying $\sup b = \int a = 1$ and $\|a \ast (b \phi^j) - \phi^j\|_p < \epsilon$ for all $1 \leq j \leq J$.
We define for each $j \leq J$
$$
\pi_n^jf(t) := a \ast_s (e^{i x_n^j \cdot \gamma(s)} b(s) f(s))(t).
$$
The weak limit condition (ii), compactness of the supports of $a$ and $b$, and an application of the Dominated Convergence Theorem imply that
$$
\lim_n \|\pi_n^j f_n - a \ast (b \phi^j)\|_p = 0.
$$
Letting $\epsilon$ to $0$, it suffices to prove that
\begin{equation} \label{E:limit bound Pn}
\lim_n \|P_n\|_{L^p \rightarrow l^{\tilde p}(L^p)} \leq 1\,\,\,\,\,\,\text{for all}\,\,\,1 \leq p \leq \infty
\end{equation}
where $\tilde p = \max(p, p')$ and $P_n := (\pi_n^j)_{j = 1}^J$.
Validity of \eqref{E:limit bound Pn} is elementary for $p = 1, \infty$. By complex interpolation it suffices to prove \eqref{E:limit bound Pn} for $p = 2$, which is equivalent to proving $\lim_n \|P_n^*\|_{l^{2}(L^2) \rightarrow L^2} \leq 1$. Note that
$P_n^*F = \sum_{j = 1}^J (\pi_n^j)^*F_j$ for $F = \{F_j\}\in l^2(L^2)$ and $(\pi_n^j)^* F_j(t) = b(t) e^{i x_n^j \cdot \gamma(t)} (a*F_j)(t)$. Thus it suffices to show that
$$
\|\sum_{j = 1}^J (\pi_n^j)^*F_j\|_{L^2}^2 \leq (1 + o_n^J(1)) \sum_{j = 1}^J \|F_j\|_{L^2}^2
$$
(where, of course, $o_n^J(1)$ is independent of $F$ and $\lim_{n \rightarrow \infty} o_n^J(1) = 0$).  
Now
$$
\|\sum_{j = 1}^J(\pi_n^j)^* F_j\|_{L^2}^2 \leq \sum_{j = 1}^J \|F_j\|_{L^2}^2 + 2 \sum_{j \neq j'} |\langle \pi_n^{j'} (\pi_n^j)^*F_j, F_{j'} \rangle|,
$$
and thus it is enough to prove that $\pi_n^{j'} (\pi_n^j)^* \rightarrow 0$ in the operator norm topology. We have $\pi_n^{j'} (\pi_n^j)^* g(u) = \int K_n^{j j'}(u, s) g(s) ds$ where
$$
K_n^{j j'}(u, s) = \int a(u - t) e^{i (x_n^j - x_n^{j'}) \cdot \gamma(t)} [b(t)]^2 a(t - s) dt.
$$
Since both $a$ and $b$ have compact support and $\supp_u K_n^{j j'}(u, s) \subset \supp a + \supp b$, by stationary phase, 
$$
\lim_n \, \sup_s \|K_n^{j j'}(u, s)\|_{L_u^1} = 0.
$$
By identical reasoning 
$$
\lim_n \, \sup_u \|K_n^{j j'}(u, s)\|_{L_s^1} = 0.
$$
Therefore $\pi_n^{j'}(\pi_n^j)^*: g \mapsto \int K_n^{j j'}(u, s) g(s) ds$ goes to $0$ in the operator norm topology. This completes the proof.
\end{proof}

\section{$L^p$ convergence}
We are finally ready to prove Theorem~\ref{mainthm}. Let $\{f_n\}$ be an extremizing sequence and $\eps > 0$. By Proposition~\ref{P:freq loc} after applying an appropriate sequence of symmetries, there exists $R = R_\eps$ such that for all sufficiently large $n$
$$
\|\scriptE f_n^R\|_q \geq B_p - \eps.
$$
Restricting $R$ to lie in the positive integers, we may apply Proposition~\ref{P:profile} along a subsequence (which is independent of $R$) to decompose $f_n^R = \sum_{j = 1}^J e^{i x_n^R \cdot \gamma(t)} \phi^{j R} + r_n^{J R}$, $J < \infty$. Furthermore since $\{f_n^R\}$ is nearly extremizing, for each $R$, there exists some large profile $\phi^{j R}$. Indeed, for all large $n$
\begin{align}
B_p^q - 2 \eps &\notag
\leq \|\scriptE f_n^R\|_q^q - \eps
\leq\sum_{j = 1}^{\infty} \|\scriptE\phi^{j R}\|_q^q
\leq B_p^{\tilde p} \bigl(\sum_{j = 1}^{\infty} \|\phi^{j R}\|_p^{\tilde p}\bigr) \max_j \|\scriptE \phi^{j R}\|_q^{q - \tilde p}
\\
&\leq B_p^{\tilde p} \max_j \|\scriptE \phi^{j R}\|_q^{q - \tilde p} \leq B_p^q \max_j \|\phi^{j R}\|_p^{q - \tilde p}.
\end{align}
Denoting this large $\phi^{j R}$ by $\phi^R$, by Proposition~\ref{P:profile} there exists $\{x_n^R\} \subseteq \Rd$ such that $\{e^{i x_n^R \cdot \gamma} f_n^R\}$ converges weakly in $L^p$ to some function $\phi^R$. Since $\lim_R \|\phi^R\|_p = 1 = \lim_n \|f_n\|_p$, by strict convexity (see Theorem $2.5$ and the proof of Theorem $2.11$ in~\cite{LL}) we have 
$$
\lim_{R} \lim_n \|f_n^R - e^{- i x_n^R \cdot \gamma} \phi^R\|_p = 0.
$$
By the triangle inequality and Proposition~\ref{P:freq loc}
\begin{equation} \label{E: fn phi R Lp}
\lim_{R} \lim_n \|f_n - e^{- i x_n^R \cdot \gamma} \phi^R\|_p = 0.
\end{equation}
Therefore for all sufficiently large $R_1, R_2$
$$
\limsup_n \|e^{-i x_n^{R_1} \cdot \gamma} \phi^{R_1} - e^{- i x_n^{R_2} \cdot \gamma} \phi^{R_2}\|_p = o_{\min(R_1, R_2)}(1). 
$$
If $\{|x_n^{R_1} - x_n^{R_2}|\}$ was unbounded for some $R_1, R_2$, then after passing through a subsequence, multiplication by $e^{i (x_n^{R_2} - x_n^{R_1}) \cdot \gamma}$ tends to zero in the weak operator topology. Thus, by H\"older's inequality,
\begin{align*}
1 &\lesssim \|\phi^{R_2}\|_p = \lim_n \tfrac{|\int(e^{i (x_n^{R_2} - x_n^{R_1}) \cdot \gamma} \phi^{R_1} - \phi^{R_2})\overline{\phi^{R_2}}| \phi^{R_2}|^{p - 2} dt|}{\|\phi^{R_2}\|_p^{p - 1}}
\\
&\lesssim \limsup_n \|e^{- i x_n^{R_1} \cdot \gamma} \phi^{R_1} - e^{- i x_n^{R_2} \cdot \gamma} \phi^{R_2}\|_p,
\end{align*}
which contradicts \eqref{E: fn phi R Lp}. Thus, for all sufficiently large $R_1, R_2$, $\{|x_n^{R_1} - x_n^{R_2}|\}$ remains bounded as $n$ goes to infinity. Applying an appropriate sequence of modulations to the $\{f_n\}$, we may assume that $\{x_n^R\}$ is bounded for all $R$. After passing to a subsequence, each $\{x_n^R\}$ converges to some $x^R \in \R$ for every sufficiently large integer $R$. Replacing $\phi^R$ with $e^{- i x^R \cdot  \gamma} \phi^R$, we may assume that $x^R = 0$ for all $R$. By \eqref{E: fn phi R Lp}, $\lim_R \limsup_n \|f_n - \phi^R\|_p = 0$, so $\{\phi^R\}$ and $\{f_n\}$ both converge in $L^p$ to some $\phi$, which must be an extremizer. This completes the proof of Theorem~\ref{mainthm}.


\bibliographystyle{plain}
\bibliography{res_moment}


\end{document}